\documentclass[10pt,leqno,twoside]{amsart}

\usepackage{amsmath,amssymb,amsfonts,amsthm}
\usepackage{fancyhdr}
\usepackage{amssymb,amsbsy}

\newtheorem{theorem}{\indent Theorem}[section]

\newtheorem{lemma}{\indent Lemma}[section]

\newtheorem*{xrem}{\indent Remark}

\newcommand\Var{\operatorname{\textrm{Var}}}

\title[Random 3-interactions]{A random graph model based on 3-interactions}

\author{\'Agnes Backhausz}
\address{Department of Probability Theory and Statistics\\Faculty of Science\\
 E\"otv\"os Lor\'and University\\P\'azm\'any P.~s.\ 1/C, H-1117
 Budapest, Hungary} 
\email{agnes@cs.elte.hu}
\author{Tam\'as F.~M\'ori}
\address{Department of Probability Theory and Statistics\\Faculty of Science\\
E\"otv\"os Lor\'and University\\P\'azm\'any P.~s. 1/C, 
H-1117 Budapest, Hungary}
\email{moritamas@ludens.elte.hu}
\dedicatory{\upshape Department of Probability Theory and Statistics,
Faculty of Science\\ E\"otv\"os Lor\'and University\\ P\'azm\'any P.~s. 1/C,
H-1117 Budapest, Hungary\\                         
\textit{E-mail address:} \texttt{agnes@cs.elte.hu, moritamas@ludens.elte.hu} \\
}
\date{1 December 2011}

\thanks{The European Union and the European Social Fund have
provided financial support to the project under the grant agreement no.\ 
T\'AMOP 4.2.1./B-09/KMR-2010-0003.}

\begin{document}

\maketitle
\begin{abstract}We consider a random graph model evolving in discrete
  time-steps that is based on $3$-interactions among
  vertices. Triangles, edges and vertices have different weights;
  objects with larger weight are more likely to participate in future
  interactions. We prove the scale free property of the model by
  exploring the asymptotic behaviour of the weight distribution.
  We also find the asympotics of the weight of a fixed vertex.\end{abstract} 

\section{Introduction}

We consider a random graph model evolving in discrete
  time-steps. The most important feature of the model is the presence
  of $3$-interactions among vertices. In our graph process vertices,
  edges and triangles get nonnegative, integer valued random weights
  which grow with time. The dynamics is driven by these weights. 

The model we are going to deal with resembles those in \cite{pub}, 
\cite{cf}, and \cite{[Mo09]}, but there are essential differences.
In \cite{cf} there is no interaction between more than two vertices,
and the weight of a vertex is simply equal to its degree. In
\cite{[Mo09]} and \cite{pub} interactions among groups of vertices do
appear, but with completely different dynamics.

Our goal is to prove that the ratio of vertices of weight $w$ tends 
to some positive constant $x_w$ almost surely, as the number of 
steps goes to infinity. We will give a recursion for $x_w$, from where 
it will be easy to see the polynomial decay of $x_w$ as 
$w\rightarrow \infty$. This is the so-called scale free property
\cite{BA}. We will also determine the asymptotics of the weight
of any fixed vertex. In the proofs martingale methods from
\cite{[Mo09]} are used. 

\section{The model}
We start with a single triangle. This has initial weight $1$, and all
its three edges have weight $1$. Later on, we will add vertices and
edges to the graph randomly. Vertices, edges, and triangles will
have nonnegative integer-valued weights, which increase according
to the random evolution of the graph.  

The graph is evolving in discrete time-steps. The sum of the weights
of triangles will be increased by $1$ at each step, while the total
weight of edges will be increased by $3$ step by step.  

At each step either a new vertex is added, which then interacts with
two already existing vertices, or $3$ old vertices interact. This
has to be decided at the beginning of the step, independently of the
past. The probability that a new vertex is born is $p$ at every step;
this is a parameter of the model. We will need $0<p\leq 1$. 

Assume that in the $n$th step a new vertex is added to the graph. We
choose two of the old vertices randomly; they will interact with
the new vertex. With probability $r$, independently of the past, the
choice is done according to the ``preferential attachment'' rule, and
with probability $1-r$ it is done ``uniformly''. $r$ is a fixed
parameter of the model. More precisely, in the case of ``preferential
attachment'' we choose the endpoints of an already existing edge
having weight $w$ with probability $\frac{w}{3n}$. Note that the sum of the
edge weights is equal to $3n$ at this moment. In the case of ``uniform''
selection two vertices are chosen with each pair having equal
probability to be selected; that is, we perform sampling without
replacement. This allows us to generate edges between old vertices
that were not connected before. 

Then the new vertex interacts with the two selected vertices. 
This means that the triangle they form comes to existence with initial
weight $1$. The two new edges connecting the new vertex to the other
two get weight $1$ each. We connect the old vertices if they are not
connected yet with an edge of weight $1$. This may only happen with
uniform selection. If the two old vertices are already connected, then
we increase the weight of that edge by $1$. To put it in another way,
we increase the weights of all three edges of the $3$-interaction by
$1$. This is the end of the step where a new vertex is generated.   

With probability $1-r$, $3$ of the old vertices will interact. With
probability $q$ they will be chosen according to the ``preferential
attachment'' rule, and with probability $(1-q)$ they will
be chosen ``uniformly''. $q$ is the third parameter of the model. This
choice is also independent of the past.  

In the ``preferential attachment'' case we choose an already existing
triangle of weight $w$ with probability proportional to its weight,
that is, with probability $\frac{w}{n}$. Note that there may exist
triangles with zero weight in the graph. A triangle has positive
weight if and only if it has already appeared in a $3$-interaction
before.  

On the other hand, in the ``uniform'' case three distinct vertices
are chosen such that every triplet has the same probability to be
selected. This is again sampling without replacement from all the
existing vertices.  

In both cases, having selected the three vertices to interact,
we draw the edges of the triangle that are not present yet. Then
the weight of the triangle is increased by $1$, as well as the weights
of the three sides of the triangle. That is, the initial weight of a
newly generated edge is $1$, while the old ones' weights are increased
by $1$.    

Now we define the weights of vertices. The weight of a vertex is the
sum of the weights of the triangles that contain it. Note that this is
just the half of the sum of weights of edges from it, because whenever
a vertex takes part in a $3$-interaction, the first sum is increased
by $1$, and the latter one is increased by $2$.   

Denote by $\mathcal F_n$ the $\sigma$-field generated by the first $n$
steps, and by $V_n$ the number of vertices after the $n$th step. Thus
$V_0=3$. Since we decide independently at each step whether a new
vertex is born, by the strong law of large numbers we obtain that
\begin{equation}\label{hae1}
V_n=pn+o\left(n^{1/2+\varepsilon}\right)\ \text{a.s.}
\end{equation}
for all $\varepsilon>0$.

Throughout this paper, for two sequences $\left( a_n\right), \left(
b_n\right)$ of nonnegative numbers, $a_n\sim b_n$ means that $b_n>0$
except finitely many terms, and $a_n/b_n\rightarrow 1$ as
$n\rightarrow \infty$.

\section{Asymptotic weight distribution}

We are interested in the distribution of weights of vertices. As we
mentioned before, this is the half of the degree of a vertex counted
with multiplicity. 

Scale-free property often emerges in models where preferential
attachment rules are applied. Therefore, throughout the paper we
suppose that the parameters do not exclude preferential attachment;
that is, either $r>0$, or $q>0$ and $p<1$.  

$X\left[ n, w\right]$ denotes the number of vertices of weight $w$
after $n$ steps. Our goal is to examine the asymptotic behaviour of
$\frac{X\left[ n, w\right]}{V_n}$; more precisely, to prove that the
ratio of vertices of weight $w$ is convergent almost surely.  
The limits are deterministic constants, which form a polynomially
decaying sequence as $w\rightarrow\infty$. We may refer to this fact as the
scale free property of the model, following the terminology of Albert
and Barab\'asi \cite{BA}. We also compute the characteristic exponent.

\begin{theorem}\label{hat1}
For $w=1,\,2,\,\dots$ we have
\[
\frac{X[n,w]}{V_n}\rightarrow x_w
\]
almost surely, as $n\rightarrow \infty$. The limits $x_w$ are
positive constants satisfying the following recursion.
\[
x_1=\frac{1}{\alpha+\beta+1}, \quad 
x_w=\frac{\alpha(w-1)+\beta}{\alpha w+\beta+1}\, x_{w-1}, 
\quad w\geq 2,
\] 
where 
\[
\alpha=\frac23\,pr+\left(1-p\right)q>0,\quad
\beta=\frac{1}{p}\bigl[2p(1-r)+3(1-p)(1-q)\bigr].
\]  
Moreover,
\[
x_w\sim C w^{-\left(1+\tfrac{1}{\alpha}\right)},
\]
as $w\rightarrow \infty$, with some positive constant $C$.
\end{theorem}

Finally we remark that $(x_d)$ is a probability distribution, its sum is 
equal to $1$.  

\proof 
First we compute the probability that a given vertex of actual
weight $w$ takes part in the 3-interaction at the $n$th step.   

If at the $n$th step a new vertex is generated and we follow the
``preferential attachment'' rule, then this probability is equal to
$\frac{2w}{3n}$, since the total weight of edges is $3n$, and the sum
of weights of edges from the given vertex is just the double of its
weight $w$. 

On the other hand, there are $\binom{V_{n-1}}{2}$
pairs of vertices, and every vertex is contained in $V_{n-1}-1$ of
them. Hence the probability of being chosen at uniform selection is
$\frac{2}{V_{n-1}}$.  

Now let us examine the case when old vertices interact. A vertex of
weight $w$ is contained in triangles of total weight $w$ by
definition, while the total sum of triangles is equal to $n$ after
$n-1$ steps. Hence the probability of being chosen is $\frac{w}{n}$
by the ``preferential attachment'' rule. With ``uniform selection'' it
is clearly $\binom{V_{n-1}-1}{2}/\binom{V_{n-1}}{3}=\frac{3}{V_{n-1}}$. 

Putting this together we get that the probability that a vertex of
weight $w$ takes part in the interaction of step $n$ is given by the
following. 
\begin{equation}\label{hae1.0}
p\left[r\,\frac{2w}{3n}+(1-r)\,\frac{2}{V_{n-1}}\right]+
(1-p)\left[q\,\frac wn+(1-q)\,\frac{3}{V_{n-1}}\right]= 
\frac{\alpha w}{n}+\frac{\beta p}{V_{n-1}}\,.
\end{equation} 

Now we determine the conditional expectation of $X[n,w]$
with respect to $\mathcal F_{n-1}$. The weights can change at most by
$1$. After $n-1$ steps we have $X\left[ n-1, w\right]$ vertices of
weight $w$. Each of them increases its weight with the probability
given above; while vertices of weight $w-1$ will count if they take
part in the 3-interaction at step $n$. Using the additive property of
expectation we obtain that  for $n\geq 1,\, w\geq 1$ the following
holds. 
\begin{multline}\label{hae1.1}
\mathbb E(X[n,w]\,|\,\mathcal F_{n-1})=
X[n-1,w]-X[n-1,w]\left[\frac{\alpha w}{n}+\frac{\beta p}{V_{n-1}}\right]\\
\hspace{2.8cm}
+X[n-1,w-1]\left[\frac{\alpha(w-1)}{n}+\frac{\beta
p}{V_{n-1}}\right]+p\delta_{w,1}
\\=X[n-1,w]\left[1-\frac{\alpha w}{n}-\frac{\beta p}{V_{n-1}}\right]
+X[n-1,w-1]\left[\frac{\alpha(w-1)}{n}+\frac{\beta p}{V_{n-1}}\right]
+p\delta_{w,1},
\end{multline}
where $X[0,w-1]$ is meant to be zero. The last term is only present
for $w=1$, because the weight of the new vertex is $1$. 

We define the following normalizing constants.
\[
c[n,w]=\prod_{i=1}^{n-1} \left(1-\frac{\alpha w}{i}-
\frac{\beta p}{V_{i-1}}\right)^{\!-1},\quad n\geq 1,\, w\geq 1.
\]

By equation \eqref{hae1}, with any positive $\varepsilon$ less than
$\frac12$ we have 
\begin{multline*}
\log c[n,w]=\sum_{i=1}^{n-1} -\log \Biggl(1-\frac{\alpha w}{i}-
\frac{\beta}{i+o\left(i^{1/2+\varepsilon}\right)}\Biggr)\\
=\sum_{i=1}^{n-1}\left(\frac{\alpha w}{i}+\frac{\beta}{i}
+o\bigl(i^{-3/2+\varepsilon}\bigr)\right)=
(\alpha w+\beta)\sum_{i=1}^{n-1}\frac{1}{i}+O(1)
\end{multline*} 
almost surely, where the error term converges as
$n\rightarrow\infty$. This implies that  
\begin{equation}\label{hae2}
c[n,w]\sim a_w n^{\alpha w+\beta}\ \text{ a.s.}
\end{equation} 
as $n\rightarrow \infty$, where $a_w$ is a positive random variable. 

Introduce $Z[n,w]=c[n,w]X[n,w]$, $n\geq 1,\, w\geq 1$. From equation
\eqref{hae1.1} it is clear that $\bigl(Z[n,w],\,\mathcal F_n\bigr)$ is
a nonnegative submartingale for every positive integer $w$. Consider
the Doob--Meyer decomposition $Z[n,w]=M[n,w]+A[n,w]$, where $M[n,w]$ 
is a martingale and $A[n,w]$ is a predictable increasing
process. Based on equation \eqref{hae1.1} we have
\begin{multline}\label{hae2.1}
A[n,w]=\mathbb E Z[1,w]+\sum_{i=2}^n \bigl(\mathbb E(Z[i,w]\,|\,
\mathcal F_{i-1})-Z[i-1,w]\bigr)\\
=\mathbb E Z[1,w]+\sum_{i=2}^n  c[i,w]\left( X[i-1,w-1]\left(
\frac{\alpha(w-1)}{i}+\frac{\beta p}{V_{i-1}}\right)+p\delta_{w,1}\right).
\end{multline}

We will also need a bound on the variation of the martingale part. By
using equation \eqref{hae2} we obtain that   
\begin{multline}\label{hae2.2}
B[n,w]=\sum_{i=2}^n \Var(Z[i,w]\,|\,\mathcal F_{i-1})=
\sum_{i=2}^n c[i,w]^2\Var(X[i,w]\,|\,\mathcal F_{i-1})\\
=\sum_{i=2}^n c[i,w]^2\Var\bigl(X[i,w]-X[i-1,w]\bigm|\mathcal
F_{i-1}\bigr)\\ 
\leq \sum_{i=2}^n c[i,w]^2\, \mathbb E \left(\left.\bigl(X[i,w]-
X[i-1,w]\bigr)^2\right\vert \mathcal F_{i-1}\right)\\
\leq 9 \sum_{i=2}^n c\left[i, w\right]^2=
O\left( n^{2(\alpha w+\beta)+1}\right).
\end{multline}
First we used the facts that $c[i,w]$ is measurable with respect to
$\mathcal F_{n-1}$, and, since there is exactly one $3$-interaction at 
each step, the change $X$ of is less than or equal to $3$. Note that
$B[n,w]$ is just the increasing process in the Doob--Meyer
decomposition of $M[n,w]^2$.

The proof continues by induction on $w$. For $w=1$ we obtain that   
\begin{equation}\label{hae3}
A[n,1]\sim p\sum_{i=2}^n c[i,1]\sim p\sum_{i=2}^n a_1 i^{\alpha+\beta}
\sim p\cdot\frac{a_1}{\alpha+\beta+1}\cdot n^{\alpha+\beta+1}
\end{equation}
almost surely, as $n\rightarrow\infty$.

On the other hand, $B[n,1]=O\left( n^{2(\alpha+\beta)+1}\right)$, 
hence, by applying Proposition VII-2-4 of Neveu \cite{[Ne75]} to
$M[n,w]$ we get that
\[
M[n,w]=o\left(B[n,1]^{1/2}\log B[n,1]\right)=
o\bigl(A[n,1]\bigr)
\]
(see Section 6 of \cite{[Mo09]} for more details of this
argument). Finally we obtain that   
\[
Z[n,1]\sim A[n,1]\ \text{ a.s.}
\]
as $n\rightarrow\infty$. Using the asymptotics of $c\left[n,1\right]$
and $A\left[n,1\right]$, that is, equations \eqref{hae2} and
\eqref{hae3}, then equation \eqref{hae1} and the definition of
$Z\left[n,1\right]$, we get that  
\[
\frac{X[n,1]}{V_n}=\frac{Z[n,1]}{c[n,1]V_n}\sim 
\frac{\dfrac{a_1}{\alpha+\beta+1}\,pn^{\alpha+\beta+1}}
{a_1 n^{\alpha+\beta}\,pn}\rightarrow
\frac{1}{\alpha+\beta+1}
\]
almost surely, as $n\rightarrow\infty$. 

Thus the theorem holds for $w=1$ with $x_1=\frac{1}{\alpha+\beta+1}$.

Suppose that the statement of Theorem \ref{hat1} holds for $w-1$ for
some fixed $w\geq 2$; that is, the ratio of vertices of weight $w-1$
converges to some constant $x_{w-1}$. By using this fact we can
compute the asymptotics of $A[n,w]$. From \eqref{hae2.1} we have 
\begin{equation*}
\begin{aligned}
A[n,w]&\sim\sum_{i=2}^n c[i,w]X[i-1,w-1]\left(\frac{\alpha(w-1)}{i}+
\frac{\beta p}{V_{i-1}}\right)\\
&\sim\sum_{i=2}^n a_w\,i^{\alpha w+\beta}\,x_{w-1}\,
V_{i-1}\left(\frac{\alpha (w-1)}{i}+\frac{\beta p}{V_{i-1}}\right)\\
&\sim a_w x_{w-1} \sum_{i=2}^n p \left(\alpha(w-1)+\beta\right) 
i^{\alpha w+\beta}\\
&\sim \frac{a_w x_{w-1} p \left (\alpha(w-1)+\beta\right)}
{\alpha w+\beta+1}\,n^{\alpha w+\beta+1} 
\end{aligned}
\end{equation*}
almost surely, as $n\rightarrow\infty$. Here we also used that
$\alpha$ and $\beta$ are both nonnegative, which is clear from their
definition. 

From inequality \eqref{hae2.2} we know that $B[n,w]=O\left( 
n^{2(\alpha w+\beta)+1}\right)$, thus Proposition VII-2-4 of 
\cite{[Ne75]} can be applied again to conclude that 
\[
M[n,w]=o\left(B[n,w]^{1/2}\log B[n,w]\right)=
o\bigl(A[n,w]\bigr).
\]

We end up with
\[
X[n,w]\sim \frac{\dfrac{a_w
x_{w-1}p(\alpha(w-1)+\beta)}{\alpha w+\beta+1}\,n^{\alpha w+\beta+1}}
{a_w\,n^{\alpha w+\beta}}
= x_{w-1}\,\frac{\alpha(w-1)+\beta}{\alpha w+\beta+1}\,np
\]
almost surely, as $n\to\infty$. Hence
\[
\lim_{n\to\infty}\frac{X[n,w]}{V_n}=
x_{w-1}\,\frac{\alpha(w-1)+\beta}{\alpha w+\beta+1}\quad\text{a.s.}
\]
Thus the induction step is complete: $x_w$ exists, and it is positive and
finite.   

Furthermore, we have a recursion for $x_w$, from where 
\begin{multline}\label{explicite}
x_w=x_1\prod_{j=2}^w \frac{\alpha(j-1)+\beta}{\alpha j+\beta+1}=
\frac{1}{\alpha w+\beta+1}\prod_{j=1}^{w-1}\frac{j+\tfrac{\beta}{\alpha}}
{j+\tfrac{\beta+1}{\alpha}}\\ 
=\frac{\Gamma\left(1+\tfrac{\beta+1}{\alpha}\right)
\Gamma\left(w+\tfrac{\beta}{\alpha}\right)}
{\alpha\,\Gamma\left(1+\tfrac{\beta}{\alpha}\right)
\Gamma\left(w+\tfrac{\beta+1}{\alpha}+1\right)}\sim 
Cw^{-\left(1+\tfrac{1}{\alpha}\right)}
\end{multline}
with some positive constant $C$, as $w$ tends to infinity. 
The proof of the theorem is complete. \qed

\begin{xrem}
From \eqref{explicite} it follows that
\[
x_w=\frac{\prod_{j=1}^{w-1}(\alpha j+\beta)}
{\prod_{j=1}^{w}(\alpha j+\beta+1)}=y_w-y_{w-1},
\]
where
\[
y_w=\prod_{j=1}^w\frac{\alpha j+\beta}{\alpha j+\beta+1}\to 0,
\]
hence 
\[
\sum_{w=1}^\infty x_w=y_0=1.
\]
\end{xrem}

\section{The weight of a fixed vertex}

In this section our goal is to determine the asymptotics of the
weight of a fixed vertex. Since the weight of a vertex is just the half
of its degree when edges are counted with multiplicity, we could
reformulate our result to obtain the asymptotics of the degree. 
   
It is clear that the weights of the vertices of the starting triangle
are interchangeable, therefore it is not necessary to deal with all
the three. Let only one of them be labelled by $0$, the other two will
remain unlabelled. Moreover, let the further vertices get labels $1$,
$2$, etc, in the order they are added to the graph. Let
$W\left[n,j\right]$ be the weight of vertex $j$ after step $n$, 
provided it exists. Otherwise let $W[n,j]$ be equal to
zero. Obviously, vertex $j$ cannot exist before step $j$. Let $I[n,j]$
denote the indicator of the event $\{W[n,j]>1\}$. 

Let us introduce the sequences
\[
b_n=\prod_{i=1}^{n}\left(1+\frac{\alpha}{i}\right)^{\!-1},\quad 
d_n=\beta p\sum_{i=1}^{n}\frac{b_{i}}{V_{i-1}}\,,
\]
with $\alpha$, $\beta$ defined in Theorem \ref{hat1}. Note that $b_n$
is deterministic, while $d_n$ is random, but $\mathcal
F_{n-1}$-measurable.   

\begin{lemma}\label{hal1}
Let $j$ and $k$ be fixed integers, $0\le j\le k$, and let
$Z[n,j]=b_nW[n,j]-d_n$. Then $\bigl(Z[n,j]I[k,j],\,\mathcal 
F_n\bigr)$ is a martingale for $n\ge k$.
\end{lemma}
\proof 
According to equation \eqref{hae1.0}, the probability that it gets new
edges at step $n+1$ is equal to  
$\dfrac{\alpha W[n,j]}{n+1}+\dfrac{\beta p}{V_n}$, provided that vertex
$j$ already exists, which surely holds if the indicator $I[k,j]$
differs from $0$. This implies that
\begin{multline*}
\mathbb E\bigl(I[k,j]W[n+1,j]\bigm|\mathcal F_n\bigr)
=I[k,j]W[n,j]+I[k,j]\left(\frac{\alpha W[n,j]}{n+1}+
\frac{\beta p}{V_n}\right)\\
=I[k,j]W[n,j]\left(1+\frac{\alpha}{n+1}\right)
+I[k,j]\frac{\beta p}{V_n}\,.
\end{multline*}
Multiplying both sides by $b_{n+1}$, we get by definition that 
\begin{multline*}
\mathbb E\bigl(b_{n+1}W[n+1,j]I[k,j]\bigm|\mathcal F_n\bigr)
=I[k,j]\left(b_n W[n,j]+b_{n+1}\frac{\beta p}{V_n}\right)\\
=I[k,j]\bigl(b_n W[n,j]-d_n+d_{n+1}\bigr),
\end{multline*}
which completes the proof of the lemma, since $d_{n+1}$ is $\mathcal
F_n$-measurable. \qed

\begin{theorem}\label{hat2} 
Fix $j\ge 0$. Then $W[n,j]\sim \zeta_j n^{\alpha}$ almost surely as 
$n\to\infty$, where $\zeta_j$ is a positive random variable. 
\end{theorem}

\proof 
First we show that this holds with a nonnegative $\zeta_j$.

Almost surely on the event that vertex $j$ exists after step $n$ we
have
\[
\mathbb P\bigl(W[n+1,j]=W[n,j]+1\bigm|\mathcal F_{n}\bigr)\ge
\frac{\alpha}{n+1}\,.
\]
Using the L\'evy-type generalization of Borel--Cantelli-lemma
\cite[VII-2-6]{[Ne75]} we get that $W[n,j]\to\infty$ with probability
$1$.  

From the definition of $b_n$ it easily follows that 
\begin{equation}\label{hae4}
b_n=\frac{\Gamma(n+1)\Gamma(1+\alpha)}{\Gamma(n+1+\alpha)}
\sim \Gamma(1+\alpha) n^{-\alpha},
\end{equation}
as $n\to\infty$.
Hence, by using equation \eqref{hae1} and the positivity of $\alpha$,
we get 
\[
d_n=\beta p\sum_{i=1}^{n}\frac{b_{i}}{V_{i-1}}=
\beta\,\Gamma(1+\alpha)\sum_{i=1}^{n-1} i^{-\alpha-1}\bigl(1+o(1)\bigr).
\]
Thus $d_n$ is almost surely convergent as $n\rightarrow
\infty$, hence the martingale of Lemma \ref{hal1} is bounded from
below. This 
martingale has bounded differences, for
\[
Z[n+1,j]-Z[n,j]\le b_n\bigl(W[n+1,j]-W[n,j]\bigr)\le b_n\le 1,
\]
and
\begin{multline*} 
Z[n,j]-Z[n+1,j]\le\bigl(b_n-b_{n+1}\bigr)W[n,j]
+\bigl(d_{n+1}-d_n\bigr)\\
\le b_{n+1}\alpha+b_{n+1}\,\frac{\beta p}{3}
\le \alpha+\frac{\beta p}{3}\,.
\end{multline*}
By Proposition VII-3-9 of \cite{[Ne75]} such a martingale either
converges or oscillates between $-\infty$ and $+\infty$, but now the
latter is excluded, hence it must converge almost surely. 

Going further, we get that $b_nW[n,j]$ is convergent almost
everywhere on the event that the weight of vertex $j$ is greater than
$1$ after step $k$. Since that weight tends to infinity,
the limit as $k\to\infty$ of this increasing sequence of events has
probability $1$. Thus $b_nW[n,j]$ is
almost surely convergent, and by equation \eqref{hae4} we get that the
statement of the theorem holds with a nonnegative $\zeta_j$.  

Now we have to prove that $\zeta_j$ is positive with probability $1$. 

In what follows, if the indicator in the numerator is zero, let us
define the fractions to be zero. Similarly to the previous lemma, for
$n\ge k$ we can write 
\begin{multline*}
\mathbb E\left(\left.\frac{I[k,j]}{W[n+1,j]-1}\,
\right\vert\mathcal F_n\right)= 
\left(\frac{\alpha W[n,j]}{n+1}+\frac{\beta p}{V_n}\right)
\frac{I[k,j]}{W[n,j]}\\+
\left[1-\left(\frac{\alpha W[n,j]}{n+1}
+\frac{\beta p}{V_n}\right)\right]\frac{I[k,j]}{W[n,j]-1}\,.
\end{multline*}
It is clear that
\[
\left(\frac{\alpha W[n,j]}{n+1}+\frac{\beta p}{V_n}\right)
\left(\frac{I[k,j]}{W[n,j]}-\frac{I[k,j]}{W[n,j]-1}\right)
\le -\frac{\alpha I[k,j]}{(n+1)\left(W[n,j]-1\right)},  
\] 
hence we get that 
\[
\mathbb E\left(\left.\frac{I[k,j]}{W[n+1,j]-1}\,
\right\vert\mathcal F_n\right)
\le\frac{I[k,j]}{W[n,j]-1}\left(1-\frac{\alpha}{n+1}\right).
\]

From this it follows that
\[
\left(\frac{e_n I[k,j]}{W[n,j]-1}\,,\ \mathcal F_n\right)
\]
is a supermartingale for $n\geq j$, where
\[
e_n=\prod_{i=1}^{n}\left(1-\frac{\alpha}{i}\right)^{\!-1}
=\frac{\Gamma(1-\alpha)\Gamma(n+1)}{\Gamma(n+1-\alpha)}
\sim \Gamma(1-\alpha)n^{\alpha}.
\]
This supermartingale is nonnegative, hence it converges almost
surely. Since 
$\lim_{k\to\infty}I[k,j]=1$ holds a.s., we obtain that 
$\dfrac{e_n}{W[n,j]-1}$ is also convergent almost surely as
$n\rightarrow \infty$. 
This implies that $\zeta_j>0$ with probability $1$, as stated.
\qed

\end{document}